\documentclass[aos,MSNbibl,nameyear,noautosecdot,dvips]{arximspdf}
\usepackage{mathbh}
\usepackage{graphicx}

\doi{10.1214/14-AOS1175REJ} %kopijuoti is PTS
\referstodoi{10.1214/13-AOS1175}
\volume{42}
\issue{2}
\pubyear{2014}
\firstpage{518}
\lastpage{531}

\makeatletter
\newcommand{\eqref}[1]{(\ref{#1})}
\def\hbeta{\hat{\beta}}
\def\tbeta{\tilde{\beta}}
\def\Exp{\operatorname{Exp}}
\makeatother

\begin{document}
\begin{frontmatter}

\title{Rejoinder: ``A significance test for the lasso''}
\runtitle{Rejoinder}

\begin{aug}
\author[A]{\fnms{Richard} \snm{Lockhart}\ead[label=e1]{lockhart@sfu.ca}},
\author[B]{\fnms{Jonathan} \snm{Taylor}\ead[label=e2]{jonathan.taylor@stanford.edu}},\\
\author[C]{\fnms{Ryan J.} \snm{Tibshirani}\ead[label=e3]{ryantibs@cmu.edu}}
\and
\author[D]{\fnms{Robert} \snm{Tibshirani}\corref{}\ead[label=e4]{tibs@stanford.edu}}
\runauthor{Lockhart, Taylor, Tibshirani and Tibshirani}
\affiliation{Simon Fraser University, Stanford University, Carnegie
Mellon University\\ and Stanford University}
\address[A]{R. Lockhart\\
Department of Statistics\\
\quad and Actuarial Science\\
Simon Fraser University\\
Burnaby, British Columbia V5A 1S6\\
Canada\\
\printead{e1}}
%adresu isvedimo komanda gale!
\address[B]{J. Taylor\\
Department of Statistics\\
Stanford University\\
Stanford, California 94305\\
USA\\
\printead{e2}\hspace*{36pt}}
\address[C]{R. J. Tibshirani\\
Departments of Statistics\\
\quad and Machine Learning\\
Carnegie Mellon University\\
229B Baker Hall\\
Pittsburgh, Pennsylvania 15213\\
USA\\
\printead{e3}}
\address[D]{R. Tibshirani\\
Department of Health, Research \& Policy\\
Department of Statistics\\
Stanford University\\
Stanford, California 94305\\
USA\\
\printead{e4}}
\pdftitle{Reply to discussions of ``A significance test for the lasso''}
\end{aug}

% HISTORY:
\received{\smonth{2} \syear{2014}}

% ABSTRACT

% KEYWORDS

\end{frontmatter}

We would like to thank the Editors and referees for their considerable
efforts that improved our paper, and all of the discussants for
their %positive
feedback, and their thoughtful and
stimulating comments. Linear models are central % is a ``bread and
%butter'' tool
in applied statistics, and inference for adaptive linear modeling
is an important active area of research.
Our paper is clearly not the last word on the subject.
Several of the discussants introduce novel proposals for
this problem; in fact, many of the discussions are interesting
``mini-papers'' on their own, and we will not attempt to reply to all
of the points that they raise. Our hope is that our paper and the
excellent accompanying discussions will serve as a helpful
resource for researchers interested in this topic.

Since the writing of our original paper, we have (with many our of
graduate students) extended the work considerably. Before responding
to the discussants, we will first summarize this new work because
it will be relevant to our responses.

\begin{itemize}
\item
As mentioned in the last section of the
paper, we have derived a ``spacing'' test of the global null
hypothesis, $\beta^*=0$, which takes the form
%
%e1 #&#
%
\begin{equation}
\label{eqspacings} \frac{1-\Phi(\lambda_1/\sigma)}{1-\Phi(\lambda
_2/\sigma)} \sim\operatorname{Unif}(0,1)
\end{equation}
for unit normed predictors, $\|X_i\|_2=1$, $i=1,\ldots, p$. As opposed
to the covariance test theory, this result is exact in finite samples,
that is, it is valid for any $n$ and $p$ (and so nonasymptotic). It
requires (essentially) only normality of the errors, and no truly
stringent assumptions about the predictor matrix $X$. In many cases,
the agreement between this test and the covariance test is very high;
details are in \citet{geomsignif} and \citet{spacing}.

\item The spacing test \eqref{eqspacings} is designed for the first
step of the lasso path. In \citet{spacing}, we generalize\vadjust{\goodbreak} it to
subsequent steps (this work is most clearly explained when we
assume no variable deletions occur along the path, i.e., when we
assume the least angle regression path, but can also be extended to
the lasso path). In addition, we study a more
general pivot that can be inverted to yield ``selection
intervals'' for coefficients of active variables at any step.

\item Similar ideas can be used to derive
$p$-values and confidence intervals for lasso active or inactive
variables at any \textit{fixed value} of $\lambda$; see
\citet{howlong}.

% \item In our theoretical results, we make the assumption that there
% are $k$ strong signals, strong enough so that with high probability
% these predictors are chosen first by LARS. A number of
% discussants correctly pointed out that this assumption is quite
% restrictive. The (nonasymptotic) spacing test proposed in
% \citet{spacing} does not make this assumption.

\item It should be noted that, in their most general form, both of the
above tests---the test at knot values of $\lambda$ in
\citet{spacing} and the test at fixed values of $\lambda$ in
\citet{howlong}---do not assume that the underlying true model is
actually sparse or even linear. For an arbitrary underlying mean
vector $\mu\in\mathbb{R}^n$, the setup allows for testing whether linear
contrasts of the mean are zero, that is, $\eta^T \mu=0$ for some
$\eta
\in\mathbb{R}^n$. Importantly, the choice of $\eta$ can be random,
that is, it can depend on the lasso active model at either a given step
or a given value of $\lambda$---in other words, both setups can be
used for \textit{post-selection inference}.

\item The question of how to use the sequential $p$-values from the
covariance test (or spacing test) is not a simple one.
As was also mentioned in the last section of our paper, in
\citet{fdrlasso}, we
propose procedures for dealing with the sequential hypothesis that
have good power properties, and have provable false discovery rate
control.
The simplest approach we call ``ForwardStop,'' which rejects for
steps $1,2,\ldots,\hat{k}_F$ where $\hat k_F={\max}\{k\dvtx
(1/k)\sum_1^k Y_i\leq\alpha\}$, and $Y_i = -\log(1-p_i)$.
\end{itemize}

We will now briefly respond to the discussants.

%s1 #&#
\section{\texorpdfstring{B\"uhlmann, Meier and van de Geer.}{B\"uhlmann, Meier and van de Geer}}\label{secbuhl}

We thank Professors B\"uhlmann, Meier and van de Geer for their
extensive and detailed discussion---they raise many interesting
points. Before addressing these, there are a few issues worth
clarifying.

\begin{itemize}
\item These authors rewrite the covariance test in what they claim
is an alternate form. Sticking to the notation in our original
paper, the quantity that they consider is
\begin{eqnarray*}
T(A,\lambda_{k+1}) &=& \bigl(\bigl\|y-X_A \tilde{
\beta}_A(\lambda_{k+1})\bigr\|_2^2 +
\lambda_{k+1}\bigl\|\tilde{\beta}_A(\lambda_{k+1})
\bigr\|_1 \bigr)/\sigma^2
\\
&&{} -\bigl(\bigl\|y-X \hbeta(\lambda_{k+1})\bigr\|_2^2 +
\lambda_{k+1}\bigl\|\hbeta(\lambda_{k+1})\bigr\|_1 \bigr)/
\sigma^2.
\end{eqnarray*}
(In\vspace*{1pt} B\"uhlmann et al. the quantities $A$ and $\lambda_{k+1}$ above are
written as $\hat{A}_{k-1}$ and $\hat{\lambda}_{k+1}$.) This is not actually
equivalent to the covariance statistic. Expanding the
above expression yields
\begin{eqnarray*}
&& 2 \bigl( \bigl\langle y, X\hbeta(\lambda_{k+1}) \bigr\rangle- \bigl
\langle y,X_A \tbeta_A(\lambda_{k+1}) \bigr
\rangle\bigr) / \sigma^2
\\
&&\quad {}+ \bigl\|X_A\tilde{\beta}_A(\lambda_{k+1})
\bigr\|_2^2 - \bigl\|X\hbeta(\lambda_{k+1})
\bigr\|_2^2 + \lambda_{k+1} \bigl(\bigl\|\tilde{
\beta}_A(\lambda_{k+1})\bigr\|_1 - \bigl\|\hbeta(
\lambda_{k+1})\bigr\|_1 \bigr),
\end{eqnarray*}
which is (two times) the covariance test statistic plus several
additional terms (the difference in squared $\ell_2$ norms of fitted
values, plus $\lambda_{k+1}$ times the difference in $\ell_1$ norms of
coefficients). These additional terms can be small,
especially if $\lambda_k$ and $\lambda_{k+1}$ are close together,
since in this case the solutions $\hbeta(\lambda_{k+1})$ and
$\tilde{\beta}_A(\lambda_{k+1})$ can themselves be close---however,
they are certainly not zero. We wonder whether this discrepancy has
affected the simulation results of B\"uhlmann et al., in Section~5 of
their discussion.

\item The authors note that the asymptotic null distributions that we
derive for the covariance test statistic require that $\mathbb{P}(B)
\rightarrow1$ as $n,p \rightarrow\infty$, for a particular event
$B$.
This event is defined slightly differently in Section~3.2, which
handles the orthogonal $X$ case, than it is in Section~4.2, which
handles the general $X$ case. Regardless, the event $B$ can be
roughly interpreted as follows: ``the lasso active model at the
given step $k$ converges to a fixed model containing the truth
(i.e., its active set contains the truly active variables, and
its active signs match those of the truly active coefficients).''

B\"uhlmann et al. comment that, to ensure that
$\mathbb{P}(B)\rightarrow1$, we assume a ``beta-min'' condition and an
``irrepresentable-type'' condition. However, this is not quite
correct. The main result of our paper, Theorem 3 in Section~4.2,
\textit{assumes that} $\mathbb{P}(B) \rightarrow1$, and uses an
irrepresentable-type condition to ensure that the conditions of the
critical Lemma 8 are met---namely, that each quantity
$M^+(j_k,s_k)$
diverges to $\infty$ quickly enough. There is no beta-min
condition employed here. If we were to have additionally assumed
a beta-min
type condition, then from this we could have shown that
$\mathbb{P}(B)\rightarrow1$. Instead, we left $\mathbb
{P}(B)\rightarrow1$ as a
direct assumption, for good reason: as described in the remarks
following Theorem 3, we believe there are weaker sufficient
conditions for $\mathbb{P}(B) \rightarrow1$
that do not require the true coefficients to be well separated from
zero---remember, for the event $B$ to hold we only need the
computed active set to contain the set of the true variables, not
equal it.

This distinction---between exact variable recovery and correct
variable screening---is an important one.
Figure~1 in the discussion by B\"uhlmann et al. shows empirical
probabilities of exact variable recovery by the lasso. It
demonstrates that, as the size $k_0$ of the true active set
increases, the
minimum absolute value of true nonzero coefficients must be quite
high in order for the lasso to recover the exact model with high
probability. But the story is quite different when we look at
variable screening; see our Figure~\ref{figscreening} below, which
replicates the simulation setup of B\"uhlmann et al., but now
records the empirical probabilities that the computed lasso model
contains the true model, after some number of steps $k \geq k_0$.
We can see that the story here is much more hopeful. For example, while
the underlying model with $k_0=10$ truly nonzero coefficients
cannot be consistently recovered after $k=10$ lasso steps, even
when beta-min is large (middle panel), this model is indeed
consistently contained in the computed lasso model after $k=20$
steps, even for very modest values of beta-min. What this means
for the covariance test, in such a setup: the asymptotic $\Exp(1)$
null distribution of the covariance statistic kicks in at some step
$k \geq k_0$, and we start to see large $p$-values. Then, by
failing to reject the null hypothesis, we correctly
screen out a sizeable proportion of truly inactive
variables.\footnote{To be fair, we are certain that B\"uhlmann et
al. are familiar with the screening properties of the lasso,
given some of these authors' own pioneering work on the subject.
Our intention here is to clarify the assumptions made in the
covariance test theory, and in particular, clarify what it means
to consider $\mathbb{P}(B) \rightarrow1$. B\"uhlmann et al. do discuss
variable screening, and remark that achieving such a property
in practice seems unrealistic, referring to their Figure~1
%and to \citet{hdscreen}
as supporting evidence. However, as explained above,
their Figure~1 examines the probability of exact model recovery,
and not screening.}

%f1 #&#
%
\begin{figure}%[b]

\includegraphics{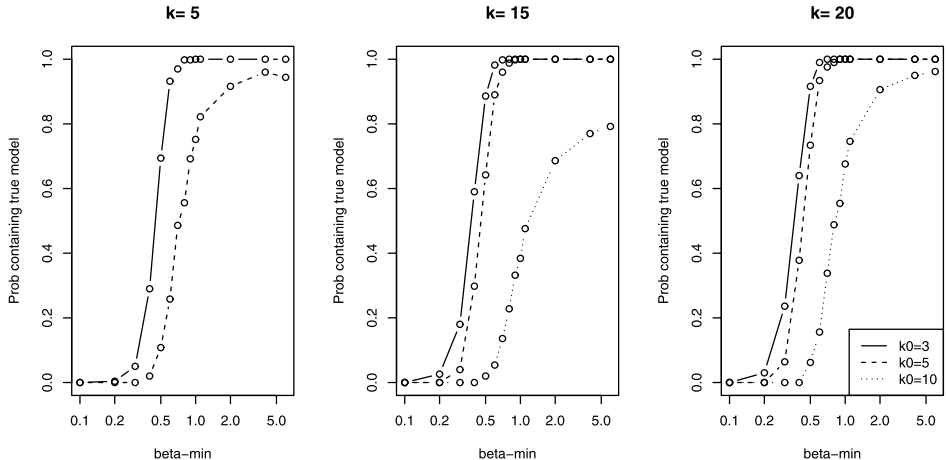}

\caption{Replication of the simulation setup considered in
the discussion by B\"uhlmann et al., but now with attention being
paid to correct variable screening, rather than exact variable
recovery. Here, $k_0$ denotes the true number of nonzero
coefficients, and $k$ the number of chosen lasso predictors
(steps along the lasso path). We see that, with high probability,
the true model is contained in the first 5, 15 or 20 chosen
predictors.}\label{figscreening}
\end{figure}

In any case, it is important to point out that the newer sequential
testing procedure in \citet{spacing} and the fixed-$\lambda$
testing procedure in \citet{howlong} do not assume a beta-min or
irrepresentable condition whatsoever, and do not require any
conditions like $\mathbb{P}(B) \rightarrow1$.
\end{itemize}

Now we respond to some of the other points raised. One of the remarks
that we made after the main result in Theorem 3 of our original paper
claims that this result can be extended to cover just the ``strong''
true variables (ones with large coefficients), and not necessarily
the ``weak'' ones (with small coefficients). B\"uhlmann et al. comment
that such an extension would likely require a ``zonal''
assumption, that bounds the number of small true coefficients, as
in \citet{hdscreen}.
As a matter of fact, we know a number of examples, with many small
nonzero coefficients, for which the conclusions of Theorem 3 continue
to hold. In any case, we emphasize that the newer sequential test of
\citet{spacing} does not need to make any assumption of this
sort, and
neither does the fixed-$\lambda$ test of \citet{howlong}.

Properly interpreting the covariance test $p$-values, as B\"uhlmann et
al. point out, can be tricky. But we believe this comes with the
territory of a conditional test for adaptive regression, since the
null hypothesis is random (and, as B\"uhlmann et al. note,
is an unobserved event). Consider the wine dataset from Section~7.1
of our original paper as an example. Looking at the $p$-values in right
panel of Table~5, one might be tempted to conclude from the $p$-value
of 0.173 in the fourth line that the constructed lasso model with
{\tt alcohol}, {\tt volatile acidity} and {\tt sulphates} contains
all of the truly active variables. There are potentially several
problems with such an interpretation (one of which being that we have
no reason to believe that the true model here is actually linear), but
we will focus on the most flagrant offense: because the constructed
lasso model is random, the $p$-value of 0.173 reflects the test of a
random null hypothesis, and so we cannot generally use it to draw
conclusions about the specific variables {\tt alcohol},
{\tt volatile acidity},
{\tt sulphates} that happened to have been selected in the current
realization. The $p$-value of 0.173 does, however, speak to the
significance of the 3-step lasso model, that is, the lasso model after 3
steps along the lasso path---said differently, we can think of this
$p$-value as reflecting the significance of the 3 ``most important''
variables as deemed by the lasso. This properly accounts for the
random nature of the hypothesis (as in any realization, the
identity of these first 3 active variables may change), and is an
example of valid post-selection inference.

Of course, one may ask: is this really what
should be tested? That is, instead of inquiring about the significance of
the 3-step lasso procedure, would a practitioner not actually want to
know about the significance of the variables {\tt alcohol},
{\tt volatile acidity} and {\tt sulphates} in particular? In a
sense, this is really a question of philosophy, and the answer
is not clear in our minds. Here, though, is a possibly helpful
observation: when we consider a single wine data set, testing the
significance of the (fixed) variables {\tt alcohol}, {\tt volatile
acidity} and {\tt sulphates} (after these 3 have been selected as
active by the lasso) seems more natural; but when we consider a
sequence of testing problems, in which we observe a new wine data set
and rerun the lasso for 3 steps on each of a sequence of days
$1,2,3,\ldots,$ testing the significance of the (random) 3-step lasso
procedure seems more appropriate.

[As an important side note, in the finite sample spacing test given in
\citet{spacing}, one can argue that both interpretations are valid,
since our inference in this work is based on conditioning on the
value of the selected variables.]

In Tables~1 and 2 of their discussion, B\"uhlmann et al. compare their
approach to
the covariance test in terms of false positive and false negative
rates. In our original paper, we had not specified a sequential
stopping rule for the covariance test, and it is not clear to us that
the two they used were reasonable. (Additionally, we are not sure
what form they assumed for the covariance test, as the representation
they present, based on the difference in lasso criterion values, is
not equivalent to the covariance statistic; see the first
clarification bullet point above.) B\"uhlmann et al. kindly sent us
their R code for their procedure, and we applied it to a subset of
their examples,
corresponding to the setup in the second row of each of their Tables~1
and 2. The results of 1000 simulations are shown in Table~\ref{tabbuhl}.
There are two setups: $n=100$, $p=80$, and $n=100$, $p=200$. In
line 1 of each, we applied their de-sparsification technique, using
the same estimate of $\sigma$ as in their
discussion. We found that this commonly overestimates
$\sigma$ by $>$100\%, so in line 2 we use the true value, $\sigma=1$.
Line 3 uses the covariance test with the ForwardStop rule of
\citet{fdrlasso}, and the true $\sigma=1$, designed to control
the FDR
at~5\%. We see that the de-spars rule does well with the inflated
estimate of $\sigma$, but produces far too many false positives when
the true value of $\sigma$ is used. Reliance on a inflated variance
estimate does not seem like a robust strategy, but perhaps there is a
way to resolve this issue. (In all fairness, B\"uhlmann et al. told us
that they are aware of this.)
The covariance test with ForwardStop does a reasonable job of
controlling the FDR, while capturing just under half of the true
signals.

%t1 #&#
%
\begin{table}
\tabcolsep=0pt
\caption{Results of a simulation study, repeating
the setup in the second row of each of Tables~1 and 2 from the
B\"uhlmann et al. discussion. Shown are the average number of
predictors called significant (out of $p=80$ or $200$), the average
number of false and true positives, the familywise error rate and
the false discovery rate}\label{tabbuhl}
\begin{tabular*}{\tablewidth}{@{\extracolsep{\fill}}@{}lccccc@{}}
\hline
& \textbf{Ave number called signif.} & \textbf{Ave FP} & \textbf{Ave TP} & \textbf{FWER} & \textbf{FDR} \\
\hline
\multicolumn{6}{@{}c@{}}{$n=100$, $p=80$} \\
(1) de-spars (estimated $\sigma$) & \phantom{0}6.89& \phantom{0}0.05& 6.84& 0.04& 0.01\\
(2) de-spars (true $\sigma$) & 17.29& \phantom{0}7.36& 9.93& 0.98& 0.43\\
(3) covTest/forwStop & \phantom{0}4.81& \phantom{0}0.25& 4.55& 0.28& 0.05
\\[3pt]
\multicolumn{6}{@{}c@{}}{$n=100$, $p=200$} \\
(1) de-spars (estimated $\sigma$) & \phantom{0}3.35 & \phantom{0}0.04& 3.30& 0.04 &0.01\\
(2) de-spars (true $\sigma$) & 44.52& 34.81& 9.71& 1.00& 0.78\\
(3) covTest/forwStop & \phantom{0}4.29 & \phantom{0}0.31& 3.97& 0.26& 0.07\\
\hline
\end{tabular*}
\end{table}

%s2 #&#
\section{\texorpdfstring{Interlude: Conditional or fixed hypothesis testing?}{Interlude: Conditional or fixed hypothesis testing}}\label{secfinal}

We would like to highlight some of the differences
between conditional and fixed hypothesis testing. This section is
motivated by the comments of B\"uhlmann et al., as well as the referees
and Editors of our original article, and personal conversations with
Larry Wasserman.

Though it has been said before, it is worth repeating:
the covariance test does not give $p$-values for classic tests of fixed
hypotheses,
such as $\beta^*_S=0$ for a fixed subset $S \subseteq\{1,\ldots, p\}$;
however, it was not designed for this purpose. As we see it:
conditional hypothesis tests like the covariance test, and fixed
hypothesis tests like that of \citet{vdgsignif} and many others (see
the references in Section~2.5 of our original paper) are two
principally different approaches for assessing significance in
high-dimensional modeling. The motivation behind the covariance test
and others is that often a practitioner becomes interested
in assessing the significance of a variable only \textit{because} it has
been entered into the active set by a fitting procedure like the lasso. If
this matches the actual workflow of the practitioner, then the
covariance test or other conditional tests seem to be best-suited to
his or her needs. A resulting complexity is that interpretation here
must be drawn out carefully (refer back to Section~\ref{secbuhl}).

On the other hand, the idea behind fixed tests like that of
\citet{zhangconf}, \citet{vdgsignif} and \citet{montahypo2},
(or at least, a typical use case in our
view) is to compute $p$-values for all fixed hypothesis $\beta_j^*=0$,
$j=1,\ldots, p$, and then perform a multiple testing correction at the
end to determine global variable significance. Even though the lasso
may have been used to construct such $p$-values, the practitioner is to
pay no attention to its output---in particular, to its active set.
And of course, the final model output by this testing procedure (which
contains the variables deemed significant) may or may not match the
lasso active set. The appeal of this approach lies in the simplicity
and transparency of its conclusions: each computed $p$-value is
associated with a familiar, classical hypothesis test, $\beta_j^*=0$
for a fixed $j$. In fact,\vspace*{2pt} we too like this approach, as it is very
direct. One drawback is that it is unclear how this might be used
for post-selection inference, if that is what is desired by the
practitioner.

We note the conditional perspective is not really a foreign one, as it is
indeed completely analogous to the (proper) interpretation of
cross-validation errors for the lasso or forward stepwise regression.
In this setting, to estimate the expected test error of a $k$-step
model computed by, say,
the lasso, we rerun the lasso for $k$ steps on a fraction of the
data set, record the observed validation error on held-out data, and
repeat this a number of times. This yields a final estimate of the
expected test error for the $k$-step lasso model; but importantly, in
each iteration of cross-validation, the selected variables will likely have
changed (since the lasso is being run on different data sets),
and so it is really only appropriate to regard cross-validation as
producing as error estimate for the $k$-step lasso procedure, not for
the particular realized model of size $k$ that was fit on the entire
data set.

Lastly, we draw attention to a connection between
our work on post-selection inference, and the de-biasing techniques
pursued by \citet{zhangconf}, \citet{vdgsignif} and
\citet{montahypo2}. In Section~7.1 of \citet{howlong}, we
show how
the framework developed in this paper can be used to form intervals or
tests for the components of a de-biased version of the true
coefficient vector, that is, something like a \textit{population
analog} of
the de-biased estimator studied by these authors. Under the
appropriate sufficient conditions [e.g., the same as those in
\citet{montahypo2}], these population de-biased coefficients converge
to the true ones, so these tests and intervals are also valid for the
underlying coefficients as well.

\section{\texorpdfstring{Buja and Brown.}{Buja and Brown}}

We thank Professors Buja and Brown for their scholarly summary
of inference in adaptive regression. We learned a great deal from it
and we enthusiastically recommend it to readers. They discuss in
detail the forward stepwise approach, and outline many
different ways to carry out inference in this setting.
To explore the $t_{\max}$ proposal
in their discussion, we carried out a simulation study.
It turns out that this is helpful in illustrating the special
properties of the covariance test with null distribution $\Exp(1/k)$,
as well as the spacing $p$-values [\citet{spacing}].

%f2 #&#
%
\begin{figure}

\includegraphics{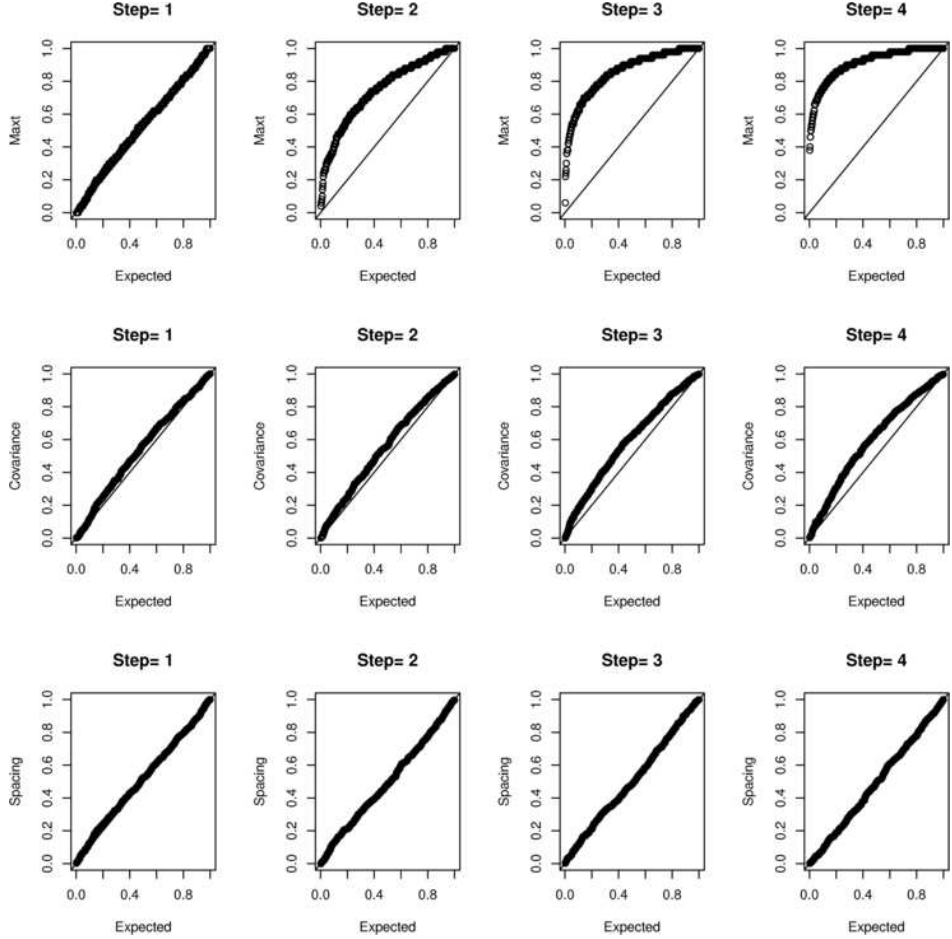}

\caption{Simulation of $p$-values for the first four
steps of using the test in \protect\eqref{eqtstat} with forward stepwise
regression (top row), the covariance test (middle row), and the
spacing test (bottom row). Details are given in the text.}\label{figbujaplot}
\end{figure}

With $n=50$, $p=10$, we generated standardized Gaussian predictors,
the population correlation between predictors $j$ and $j'$
being $0.5^{|j-j'|}$. The true coefficients were $\beta^*=0$, and
the marginal error variance was $\sigma^2=1$. The middle
and bottom panels of Figure~\ref{figbujaplot} show quantile--quantile
plots of the covariance $p$-values and spacing $p$-values for the first
four steps of the least angle regression path [see \eqref{eqspacings}
for the spacing test in the first step, and \citet{spacing} for
subsequent steps]. In the top panel, we have applied forward
stepwise regression, using the test statistic
%
%e2 #&#
%
\begin{equation}
\label{eqtstat} t_{\max}(y)=\max_{j=k,\ldots, p}
\bigl|t^{(j)}(y)\bigr| \qquad\mbox{where } t^{(j)}(y)=\frac{\langle X_{j\cdot
A},y\rangle}{\| X_{j\cdot A}\|_2},
\end{equation}
per the proposal of Buja and Brown. Here, $A$ is the set of active
variables currently in the model and $X_{j\cdot A}$ denotes the $j$th
predictor orthogonalized with respect to these variables.
Note that we have used the true value $\sigma^2=1$ in \eqref{eqtstat}
(and in the covariance and spacing tests as well). As suggested by
Buja and Brown, we simulated $\varepsilon\sim N(0,I)$ in order to
estimate the $p$-value $\mathbb{P}(t_{\max}(\varepsilon) >
t_{\max}(y))$.

All three tests look good at the first step, but the forward stepwise
test based on \eqref{eqtstat}
becomes more and more conservative for later steps.
The reason is that the covariance test and the spacing test (even
moreso) properly account for the selection events up to and including
step $k$. To give a concrete example, the forward stepwise test
ignores the fact that at the second step, the observed
$t_{\max}$ is the \textit{second} largest value of the statistic
in the data, and erroneously compares it to a null distribution of
\textit{largest} $t_{\max}$ values. This creates a conservative
bias in the $p$-value. If predictor $j$ were chosen at the first step
of the forward stepwise procedure, then a correct numerical
simulation for
$t_{\max}$ at the second step would generate
$y^*=X_j \hbeta_j + \varepsilon$ (with $\hbeta_j$ being
the least squares coefficient on variable $j$), and only keep those
$y^*$ vectors for which predictor $j$ is chosen at the first
step, using these to compute $t_{\max}(y^*)$ [which
equals $t_{\max}(\varepsilon)$]. Such a simulation setup might be
practical for a few steps, but would not be practical beyond
that, though there do exist efficient algorithms for
sampling from such distributions. Remarkably, the covariance and
spacing tests are able carry out
this conditioning analytically.

On a separate point, we agree with Buja and Brown that inferences
should not typically focus on the true regression coefficients when
predictors are highly correlated,
and even the definition of FDR seems debatable in that setting.
In \citet{uvr}, we propose an alternative definition of FDR,
called the
``Uninformative Variable Rate'' (UVR), which tries to finesse this
issue by projecting the true mean $X\beta^*$ onto the set of
predictors in the current model. A selection is deemed a false
positive if it has a zero coefficient in this projection.
For example, in a model with $\beta_1^*=5$, $\beta^*_2=0$ and
$\operatorname{Cor}(X_1,X_2)=0.95$, the selection of $X_2$ by itself would
be considered a false positive in computing the FDR. But this does not
seem reasonable, and the UVR would instead consider it a true
positive.

As Buja and Brown mentioned, we have proposed a method for
combining sequential $p$-values to achieve FDR control in
\citet{fdrlasso}. But we believe there is more to do, especially in
light of the last point just raised.
%deriving tests with good power properties.

Finally, as they remark, our tests will not be valid if the
practitioner uses them in combination with other selection techniques,
or as they put it, the data analyst is ``arbitrarily informal in their
meta-selection of variable selection methods.''
As they point out, the POSI methods they propose in \citet{posi}
are valid even in that situation. This is a very nice property, but of
course the pressing question is: are the inferences too conservative
as a result of protecting the type I error in such a broad sense?

%s4 #&#
\section{\texorpdfstring{Cai and Yuan.}{Cai and Yuan}}

We are grateful to Professors Cai and Yuan for their suggestion of an
alternative test based on the Gumbel distribution. In the most basic
setting, testing at the first step (i.e., global null hypothesis) in
the orthogonal $X$ setting, both our proposal and theirs stem from the
same basic arguments. To see this,
suppose that $V_1 \geq\cdots\geq V_p>0$ are the ordered absolute
values of a sample from a standard normal distribution. Then, as
$p \to\infty$,
%
%e3 #&#
%
\begin{equation}
\label{eqablim} b_p(V_1 - a_p) \stackrel{d} {
\rightarrow}\operatorname{Gumbel}(0,1),
\end{equation}
where
\[
a_p = \Phi^{-1} \bigl(1-1/(2p) \bigr) = \sqrt{2 \log p} -
\frac{\log\log p +\log\pi}{2\sqrt{2\log p}} + o(1/\sqrt{\log p}) %
% +\log\pi}{2\sqrt{2\log p}}+o\big(\log^{-1/2} (p)\big),
\]
and
\[
b_p = \sqrt{2\log p} \bigl(1+o(1) \bigr).
\]
We used this and the fact that
$b_p(V_1-V_2) \stackrel{d}{\rightarrow}\Exp(1)$
to handle the orthogonal $X$ case.
Dividing \eqref{eqablim} by $b_p^2$, we see that
\[
\frac{V_1+a_p}{b_p} =\frac{V_1-a_p}{b_p} + 2 + o(1) \rightarrow2
\]
and multiplying by \eqref{eqablim}, we get
\[
V_1^2-a_p^2 \stackrel{d} {
\rightarrow}\operatorname{Gumbel}(0,2),
\]
which may be rearranged to give Cai and Yuan's observation [since
$a_p^2 = 2\log p -\log\log p - \log\pi+o(1)$]. Hence for the
orthogonal case, under the global null, we are basically using the
same extreme value theory.

\begin{figure}

\includegraphics{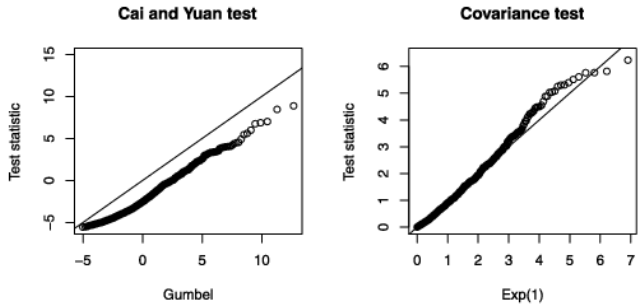}

\caption{Quantile--quantile of Gumbel test (left panel) and
covariance test (right panel) with features having pairwise
correlation $0.7$.}\label{figcai}
\end{figure}

But for a general predictor matrix $X$, even if we stick to
testing at the first
step, we believe the Gumbel test does not share the same kind of
parameter-free asymptotic behavior of the covariance test.
Specifically, take $X^TX$ to be the $p \times p$ matrix with 1 on the
diagonal and every off diagonal element equal to some fixed
$\rho\in(0,1)$. In this case, we can show that under the global null,
\[
V_1 - \sqrt{2\rho\log p}\stackrel{d} {\rightarrow}\bigl|N(0,1-\rho)\bigr|,
\]
so the asymptotic distribution depends on $\rho$,
and the procedure suggested by Cai and Yang must fail.
Figure~\ref{figcai} shows an example with $n=100$, $p=50$, and
$\rho=0.7$. The Gumbel approximation is poor, while the $\Exp(1)$
distribution for the covariance test statistic still works well.

%{\bf Richard Observes: This distribution does not depend on
% $n$. Our $X^TX$ has rank $p$ so implicitly we have $n \ge p+1$
% because we have centred the variates.}
%
%f3 #&#
%

Another important point is that, for a general $X$, the test proposed
by Cai and Yuan does not apply to the sequence of variables entered
along the lasso path. Cai and Yuan assume that,
given a current active set $A$, the variable $j$ to be entered is that
which maximizes the drop in residual sum of squares. (In their
notation, the representation $R_j = \max_{m \notin A} R_m$ is what
allows them to derive the asymptotic Gumbel null distribution for
their test.)
While this is true at each lasso step in the orthogonal $X$ case, it
is certainly not true in the general $X$ case. Meanwhile, for an
arbitrary $X$, it is true in forward stepwise regression (by
definition).

\section{\texorpdfstring{Fan and Ke.}{Fan and Ke}}

Professors Fan and Ke extend the covariance test and its null
distribution to the SCAD and MCP penalties, in the orthogonal $X$
case. This is very exciting. We wonder whether this can be extended to
arbitrary $X$, and whether the spacing test [\citet{spacing}] can be
similarly generalized.

Fan and Ke (and also B\"uhlmann, Meier and van de Geer) also study the
important issue of the power of the covariance test,
relative to the ``RSSdrop'' and ``MaxCov'' statistics.
The discussants here have honed in on
the worst case scenario for the covariance test, in which two
predictors have large and equal coefficients. In this situation, the
LARS algorithm
takes only a short step after the first predictor has been
entered, before entering the second predictor, and the hence the
$p$-value for the first step is not very small. For this reason, better
power can be achieved by constructing functions of more than one
covariance test $p$-value, as illustrated by Figure~4 in the discussion
of Fan and Ke. We note, however, that neither RSSdrop nor MaxCov have
tractable null distributions in the general $X$ case, and it is not
even clear how to approximate these null distributions by simulation
except in the global null setup. Power concerns were also part of the
motivation for our development of the sequential tests in
\citet{fdrlasso}. Also, it is worth mentioning that the framework of
\citet{spacing} actually allows for combinations of the knots
$\lambda_j$, $j=1,\ldots, k$ from the first $k$ steps, so that we could
form an exact test based on, for example, $\sum_{j=1}^k
\lambda_j$ is
this was seen to have better power. Overall, the issue of the ``most
powerful sequential test'' remains an open and important one.

Continuing on the topic of power, Fan and Ke (and again, B\"uhlmann et
al.) raise asymptotic concerns.\vspace*{1pt} They suggest that power against
coefficients on the
order $O(n^{-1/2})$ is desirable. A first clarification: if elements of
$y$ and the rows of $X$ are generated by i.i.d. sampling, then the
matrix $X^TX$ grows like $n$; our standardization, in which $X^TX$
has 1 in each diagonal entry, corresponds to multiplying $\beta^*$ by
$\sqrt{n}$ in this i.i.d. sampling context. The rate $O(n^{-1/2})$
mentioned then becomes
$O(1)$. Power results will generally depend on $X$, and a complete
discussion would be outside of the scope of this discussion, but some
insight into what is possible or what is reasonable to expect may be
gained by considering the orthogonal
case. Consider now the problem of testing the global null against the
alternative $\beta^*_{j_0} \neq0$ and $\beta^*_j=0$ for all
$j\neq j_0$, with $j_0$ known. For $|\beta^*_{j_0}|=\nu$,
fixed, we get nontrivial limiting power by rejecting if
$|U_{j_0}| = |X_{j_0}^Ty|>z_{\alpha/2}$, as usual.\vspace*{2pt} But
realistically, ${j_0}$~will
not be known and it will be sensible to ask about the average power
over all ${j_0}\in\{1,\ldots,p\}$. The problem of testing $\beta^*=0$
against the hypothesis that there is a unique ${j_0}\in\{1,\ldots,p\}
$ for
which $\beta^*_{j_0} \neq0$ is invariant under permutations
of the entries in $U$. Let ${\mathcal T}_p$ denote the class of all
permutation invariant tests $T(U)$; our test $T_1$ and any other tests
which are functions of the order statistics of $U_j$, $j=1,\ldots, p$ are
permutation invariant. Let $B_p(\nu)$ be the set of $\beta^*$ with
exactly one nonzero entry satisfying $|\beta^*_j| \le\nu$. We can
prove that if $\nu_p$ is any sequence of constants with
\[
\sqrt{2\log p }-\nu_p \to\infty, %
\]
then
\[
\sup_{T \in{\mathcal T}_p, \beta^*\in B_p(\nu_p)} \bigl|\operatorname
{Power}\bigl(T,\beta^*\bigr) -
\operatorname{Level}(T)\bigr| \to0. %
\]
For tests which are not permutation invariant, we can prove
\[
\sup_{T, \beta^*\in B_p(\nu_p)} \bigl|\operatorname{AveragePower}\bigl
(T,\beta^*\bigr)
- \operatorname{Level}(T)\bigr| \to0, %
\]
where now AveragePower denotes, for a given $\beta^*\in B_p(\nu_p)$,
the average over the $p$ vectors obtained by permuting the
entries of $\beta^*$. In other words, unless $\beta^*$ has an entry on
the order of $\sqrt{2\log p}$, there is no permutation invariant way to
distinguish the null from the alternative. On the other hand, if
$a_p=\sqrt{2\log p} - \log(\log p)/(2\sqrt{2\log p})$ and
\[
a_p(a_p-\nu_p) \to-\infty, %
\]
then our test has limiting power 1 in this context.
This $\sqrt{2\log p}$ rate, then, cannot be substantially improved in
general. The same conclusion holds if $B_p(\nu)$ is replaced by the
intersection of the $O(1)$ ball $\{\beta\dvtx\|\beta\|_2 \le\Delta\}$
with $\{\beta\dvtx|\bar\beta| \le\varepsilon_p/\sqrt{p}\}$.
Here $\Delta$ is any fixed constant,
$\bar\beta= \sum_{j=1}^p \beta_j/p$,
and $\varepsilon_p$ is any sequence shrinking to 0. Notice that if
$\beta^*$ in this set is known then using a likelihood ratio test
against that alternative achieves nontrivial asymptotic power
(provided $\|\beta^*\|_2$ stays away from 0). If the permutation group
is expanded to the signed permutation group, then the condition on
$\bar\beta$ may be deleted; natural procedures will have this added
sign invariance in the orthogonal case.

%s6 #&#
\section{\texorpdfstring{Lv and Zheng.}{Lv and Zheng}}

Professors Lv and Zheng explore extensions of these ideas to nonconvex
objective functions, for example, a combination of Lasso and the SICA
penalty. This is interesting but seems difficult, as even the
computation of the global solution is infeasible in general. However,
the existing asymptotic results for these methods suggest that
inference tools might also prove to be tractable. Regarding the
significance of each active predictor conditional on the set of all
remaining active predictors: the spacing theory in \citet{spacing}
provides a method for doing this.

%{\bf Added by RAL in a very tentative voice.}
Lv and Zheng also suggest extra shrinkage, replacing
$\lambda_{k+1}$ in our, and their, test statistics by
$c\lambda_{k+1}$, in the hopes that a better choice of $c<1$ will lead
to an improved $\Exp(1)$ approximation. In knot form, this would look
like %The knot form of
%our statistic shows that we are considering, in many cases,
%
\[
C(A,s_A,j,s)\lambda_k(\lambda_k - c
\lambda_{k+1}) = T_k +(1-c)C(A,s_A,j,s)
\lambda_k\lambda_{k+1}. %
\]
Typically, $\lambda_k$, $\lambda_{k+1}$ are drifting to $\infty$ with
$p$, so the shrinkage factor $c$ will have to be chosen carefully in
order to control the second term above;
it seems that $c \rightarrow1$ is needed whenever the limit of
$T_k$ is $\Exp(1)$.
%Whenever $C(A,s_A,j,s) = 1$, as below Lemma 1 in our
%original article, we will need $c$ quite near 1.

%s7 #&#
\section{\texorpdfstring{Wasserman.}{Wasserman}}

Professor Wasserman appropriately points out the stringency of
assumptions made in our paper, assumptions that
are common to much of
the theoretical work on high-dimensional regression. We would like to
reiterate that three of the offending assumptions in his
list---that is, the assumptions that the true model is linear and is
furthermore
sparse, and that the predictors in $X$ are weakly correlated---are not
needed in the newer works of \citet{spacing} and \citet{howlong}.
In general, though, we do agree that the rest of assumptions in his list
(implying independent, normal, homoskedastic errors) are used for as a
default starting point for theoretical analysis, but are certainly
suspect in
practice.

Wasserman outlines a model-free approach to inference in adaptive
regression based on sample-splitting and the increase in predictive
risk due to setting a coefficient to zero. The proposal is simple and
natural, and we can appreciate model-free approaches that use sample
splitting like this one. However, we worry about
the loss in power due to splitting the data in half, especially when
$n$ is small relative to $p$. As he says, this may be the
price to
pay for added robustness to model misspecification. How steep is
this price? It would be interesting to investigate, both
theoretically or empirically, the precise power lost due to sample
splitting.
Also, we note that
the random choice of splits will also influence the results, perhaps
considerably. Therefore, one would need to take multiple random
splits, and somehow combine the results at the end;\vadjust{\goodbreak} but then the
interpretation of the final ``conditional'' test seems challenging.
We are eager to read a completed manuscript on this interesting idea.

His discussion of conformal prediction is fascinating; this is an area
completely new to us. And finally, we thank him for his clearly
expressed reminder of the difficulties of determining causality from a
standard statistical model.

%s8 #&#
\section{\texorpdfstring{Thanks.}{Thanks}}

We thank all the discussants again for their contributions. They
have given us much to think about. We hope that our original paper,
the subsequent discussions and this response will be a valuable
resource for researchers interested in inference for adaptive
regression.


%suskaldyti doi

% imsref loaded by linak, 2014-03-18 16:09:40
%
% imsref loaded by linak, 2014-03-19 13:25:10

\printaddresses

\end{document}